\documentclass[twoside]{article}

\usepackage[preprint]{aistats2027}
\usepackage{amsfonts}
\usepackage{url}
\usepackage{amsthm}
\usepackage{comment}
\usepackage{graphicx}

\newtheorem{theorem}{Theorem}
\newtheorem{proposition}[theorem]{Proposition}

\newtheorem{corollary}[theorem]{Corollary}
\newtheorem{remark}{Remark}
\usepackage[numbers,sort&compress]{natbib}

\usepackage[colorlinks=true, allcolors=blue]{hyperref}
\usepackage{xcolor}
\usepackage{algorithm}
\usepackage{algpseudocode}
\usepackage{booktabs}
\usepackage{multirow}
\usepackage{subcaption}
\begin{document}

\twocolumn[
\aistatstitle{ Taming Nonconvexity with Epigraph Lifting
}
 
\aistatsauthor{ Clara Lage \And Hok Shing Wong  \And Matthias Ehrhardt   }

\aistatsaddress{
 \parbox[t]{0.33\textwidth}{\centering\small 
    Institut Camille Jordan\\ 
    \'Ecole Centrale de Lyon\\ 
    \texttt{clara.lage@ec-lyon.fr}
  }

  \And

  \parbox[t]{0.33\textwidth}{\centering\small 
    LMU Munich\\ 
    Munich Center for Machine Learning\\ 
\texttt{hok.wong@math.lmu.de}
  }

  \And

  \parbox[t]{0.33\textwidth}{\centering\small 
    University of Bath\\     \texttt{m.ehrhardt@bath.ac.uk}
  }
 }
]

\begin{abstract}
We consider inverse problems regularized by a trained neural network with \emph{convex activations} and \emph{unrestricted weights}. The resulting optimization problem is generally nonconvex, and gradient-based methods may converge to poor local minima. We introduce a \emph{penalized epigraph lifting} of the network and show that the lifted problem admits a \emph{difference-of-convex (DC) decomposition}. For Lipschitz continuous activations, we establish \emph{finite exactness}: the original objective is recovered exactly at a finite and explicitly computable penalty threshold $\bar\rho$. Building on this structure, we propose \emph{H-DC} (Homotopic DC), a homotopic scheme that solves a sequence of DC subproblems with increasing $\rho$, warm-starting each level from the previous solution. Numerical experiments on synthetic problems, a manifold denoiser regularizer, and image inpainting with a learned Mumford--Shah regularizer show that H-DC consistently escapes poor local minima and reaches solutions with lower objective values than gradient-based methods from the same initialization.
\end{abstract}

\section{Introduction}\label{sec:intro}
A growing class of problems in machine learning and signal processing requires solving an optimization problem in which
a neural network appears inside the objective with fixed weights~\cite{carlini2017towards,salman2019convex,xia2022gan,zhang2018efficient}, a setting we refer to as
\emph{network inference}. A particular instance arises in the context of learned regularization in inverse problems, where a neural network defines a learned prior and reconstruction is obtained by solving
\begin{equation}
    \min_{x} \; \left\{D(x) + \lambda R_\theta(x)\right\},
    \tag{P}
    \label{eq:P}
\end{equation}
where $D$ is a data-fidelity term and $R_\theta$ is a network-parametrized regularizer \cite{arridge2019solving, goujon2023neural}. Networks also enter regularizers through \emph{least-squares residual (LSR) regularization}, where a learned denoiser $N_\theta$
defines a coupling term $\frac{\alpha}{2}\|x -
N_\theta(x)\|^2$ \cite{zou2023deep}. In each case, once the weights are
fixed, inference reduces to a (generally nonconvex)
optimization problem over $x$.

\paragraph{Expressivity versus tractability.}
The flexibility of neural parameterizations comes with a
fundamental optimization challenge. Restricting network architectures to ensure convexity or controlled
nonconvexity leads to tractable inference but limits the class of representable priors. In pursuit of richer expressivity, we consider feedforward networks with \emph{convex activations and unrestricted weights}, a setting
which is generally nonconvex and where optimization offers limited guarantees.

\subsection{Related Work}
\label{sec:related}

\paragraph{Learned regularizers.}
Learned regularizers parameterized by neural networks have become a popular approach for solving inverse
problems~\cite{goujon2023neural,hertrich2025learning,
lunz2018adversarial,roth2009fields}. \emph{Convex} architectures such as ICNNs~\cite{amos2017input} lead
to tractable inference with global convergence
guarantees~\cite{mukherjee2020learned}. \emph{Weakly convex}
regularizers~\cite{goujon2024learning,shumaylov2024weakly} relax convexity while retaining stationarity guarantees. More complex definitions of $R$ through a pretrained denoiser $N_\theta$ arise in \emph{regularization by denoising (RED)}~\cite{romano2017little} and \emph{LSR regularization}~\cite{zou2023deep}. Our work addresses
networks with convex activations and unrestricted weights, generalizing the convex and weakly
convex settings while accommodating these more general formulations.

\paragraph{Lifting and epigraph relaxations.}
Introducing auxiliary
variables for the intermediate layer outputs and relaxing the forward-pass equalities, a technique
known as \emph{lifting}, has been extensively studied in the context of network training~\cite{askari2018lifted,carreira2014distributed,gu2020fenchel,wang2023lifted}. Recently, lifting was applied to ICNNs in~\cite{ehrhardt2026primaldual}, where the forward-pass equalities are relaxed to epigraph inequalities. The nonnegativity weight constraints ensure that the epigraph lifting is tight at any minimizer and the lifted problem is equivalent to the original without any penalization. Our work extends this to general networks via a penalty parameter that closes the tightness gap at a finite and
computable threshold.

\paragraph{DC programming and homotopy.}
\emph{DC programming}~\cite{lethi2012exact} and its boosted variant BDCA~\cite{artacho2020positive} provide
principled algorithms for difference-of-convex problems.
\emph{Homotopy methods}~\cite{iwakiri2022single,
sato2025explicit} tackle nonconvex optimization through
a sequence of easier problems, typically via external smoothing. Our continuation path is instead induced by the network architecture itself, with a finite exactness
guarantee.

\paragraph{Nonconvex regularization in imaging.}
The \emph{truncated quadratic}~\cite{mumford1989optimal,
strekalovskiy2014real} and \emph{$\ell_p$
penalties}~\cite{soubies2015continuous} are widely used
nonconvex regularizers in imaging. DCA has been applied to truncated quadratic regularization~\cite{deng2021preconditioned}, and exact relaxations of $\ell_0$ penalties are studied
in~\cite{essafri2025exact}. Our work illustrate that learned
approximations of such regularizers can be optimized
more effectively via epigraph lifting.
 
\paragraph{Epigraph lifting for network inference.}
We propose to address the network inference problem by exploiting the internal layered structure of the network directly. \emph{Lifting} the problem by introducing auxiliary variables and relaxing forward-pass equalities to \emph{epigraph inequalities} yields an exact formulation for convex networks, but breaks down for general networks: negative weights allow the lifted variables to grow arbitrarily large while remaining feasible, making the unpenalized problem degenerate. We introduce a \emph{penalty parameter} $\rho > 0$ on the epigraph slack that recovers the original objective exactly for all $\rho$ above a finite and explicitly computable threshold $\bar\rho$. This yields a family $\{G_\rho\}_{\rho \geq 0}$ satisfying $G_\rho \leq R$ for all $\rho$, with $G_\rho = R$ for $\rho \geq \bar\rho$. The second key ingredient is that this homotopy comes with additional optimization structure: under convex activations, each lifted problem admits a \emph{difference-of-convex (DC) decomposition}. Rather than merely constructing a sequence of smoother surrogate objectives, the proposed lifting produces a sequence of structured nonconvex problems addressable by DC optimization. This provides a natural combination of two complementary mechanisms: \emph{continuation} is used to navigate progressively toward the original nonconvex landscape, while the \emph{DC structure} is exploited to solve each intermediate problem globally, yielding our method
\textbf{H-DC} (Homotopic DC).

\paragraph{Contributions.}
Our main contributions are threefold:
\begin{itemize}
    \item We introduce a \emph{penalized epigraph lifting}
    of general pretrained neural regularizers that
    generates a parameterized family $G_\rho \leq R$
    and establish a finite exactness threshold $\bar\rho$
    for which $G_\rho = R$ (Proposition~1).

    \item We use this family to define a
    \emph{network-induced homotopy} from a lifted
    relaxation to the original neural-regularized
    problem, leading to a continuation strategy based
    on successively increasing $\rho$ and warm-starting
    the resulting optimization problems.

    \item We establish a \emph{DC structure} for the
    lifted problems under convex activations and exploit
    it algorithmically along the continuation path,
    yielding \textbf{H-DC}. Numerical experiments show
    that H-DC reaches solutions with lower objective
    values than direct local optimization of the
    original nonconvex neural regularizer.
\end{itemize}

\section{Penalized Epigraph Lifting}
\label{sec:lifting}

\subsection{Network inference problem}\label{sec:net-inf}

We consider a feedforward neural network defined recursively by
\begin{equation}
\begin{aligned}
    z_0 &= x, \\
    z_\ell &= \sigma_\ell\!\left(
        W_\ell z_{\ell-1} + U_\ell x + b_\ell
    \right),
    \quad \ell=1,\ldots,L,\\
    R(x) &= V^\top z_L + c,
\end{aligned}
\label{eq:network}
\end{equation}
where $W_\ell$ and
$U_\ell$ are weight matrices,
$b_\ell$ are biases, $V$
and $c$ define the output layer, and $\sigma_\ell$
are componentwise activation functions. Throughout, we
assume that each $\sigma_\ell$ is \emph{convex}, which includes
ReLU, leaky ReLU, and smoothed ReLU as special cases. In this
paper, we focus on the variational reconstruction framework~\eqref{eq:P} with
learned regularizer. While the architecture~\eqref{eq:network} generally defines a scalar-output
regularizer, it also accommodates the LSR setting through a
simple modification of the output layer.

\begin{remark}[LSR regularizer]
The least-squares residual regularizer $R(x) =
\frac{1}{2}\|x - N_\theta(x)\|^2$ fits within
architecture~\eqref{eq:network} as follows. Let $z_L = N_\theta(x)$
denote the output of the denoiser, and define an additional layer
$z_{\mathrm{LSR}} = x - z_L$. Setting $\sigma_{\mathrm{LSR}}(s) =
s^2$ componentwise and $V = \frac{1}{2}\mathbf{1}$ then gives
$R(x) = V^\top \sigma_{\mathrm{LSR}}(z_{\mathrm{LSR}}) =
\frac{1}{2}\|x - N_\theta(x)\|^2$.
\end{remark}

\subsection{Epigraph lifting}\label{sec:epi-lift}

The key observation is that, for convex $\sigma_\ell$, the
epigraph
\[
    \mathrm{epi}(\sigma_\ell) = \bigl\{
        (s, t) : t \geq \sigma_\ell(s)
    \bigr\}
\]
is a convex set. We introduce auxiliary variables $z_1, \ldots,
z_L$ for the intermediate layer outputs and relax the forward-pass
equalities $z_\ell = \sigma_\ell(W_\ell z_{\ell-1} + U_\ell x +
b_\ell)$ to epigraph inequalities, defining the feasible set
\begin{equation}
\mathcal{C}_\sigma^L := \bigl\{
    (x, \boldsymbol{z}) :
    z_\ell \geq \sigma_\ell(\phi_\ell), \ell = 1, \ldots, L
\bigr\},
\label{eq:C}
\end{equation}
where $\phi_\ell = W_\ell z_{\ell-1} + U_\ell x + b_\ell$ denotes
the pre-activation at layer $\ell$ and $\boldsymbol{z} = (z_1, \ldots, z_L)$.

Note that equality $z_\ell = \sigma_\ell(\phi_\ell)$ holds if and only if $(x, z_1, \ldots, z_L)$ lies on the boundary of $\mathcal{C}_\sigma^L$, corresponding exactly to the forward pass of the network. The epigraph inequality relaxes this to allow $z_\ell$ to lie anywhere above the graph of $\sigma_\ell$.

Since each $\sigma_\ell$ is convex, $\mathcal{C}_\sigma^L$
is a convex set, which is the key structural property underlying both the feasibility constraints and the DC decomposition in Proposition~\ref{prop:DC_convex_activation_L}. This assumption places no restriction on the weights, the non-convexity of the problem comes entirely from the weights, not the activations. For ICNNs, the non-negativity weight constraints ensure that increasing any lifted variable can only increase the subsequent layers and therefore the output $V^\top z_L$. The epigraph
constraints are thus automatically tight at any minimizer, and
the unpenalized lifting already recovers $R$ exactly \cite{ehrhardt2026primaldual}. For networks with unrestricted weights, this argument breaks down: a negative weight at layer $\ell+1$ means that increasing $z_\ell$ can decrease the output, so the lifted variables may drift arbitrarily far from the forward pass while reducing $V^\top z_L$. Although the unpenalized lifted problem always provides a lower bound for $R$, it may become degenerate; either unbounded from below or yielding a trivial lower bound. Adding a penalty $\rho_\ell \mathbf{1}^\top(z_\ell - \sigma_\ell(\phi_\ell))$ at each layer discourages this drift. This gives the \emph{penalized lifted regularizer}
\begin{equation}
G_{\boldsymbol{\rho}}(x) = \inf_{\substack{\boldsymbol{z} \\ 
(x,\boldsymbol{z}) \in \mathcal{C}_\sigma^L}}
\Bigl\{
    V^\top z_L
    + \sum_{\ell=1}^{L} \rho_\ell \,
    \mathbf{1}^\top \bigl(
        z_\ell - \sigma_\ell(\phi_\ell)
    \bigr)
\Bigr\},
\label{eq:G_rho}
\end{equation}
where $\boldsymbol{\rho} = (\rho_1, \ldots, \rho_L) \geq 0$. Another possible approach to prevent such degeneracies is to directly restrict the feasible set $\mathcal{C}_\sigma^L$, as considered in~\cite{salman2019convex}. From a practical perspective, however, the penalized formulation $G_{\boldsymbol{\rho}}$ exhibits more favorable properties for optimization, as the penalty both prevents degeneracy and provides a natural continuation path toward the original objective, as we discuss in Section~\ref{sec:theory}.

\subsection{Theoretical properties}
\label{sec:theory}
The family $\{G_{\boldsymbol{\rho}}\}$ enjoys several key
properties that underpin the H-DC algorithm.

\begin{proposition}
\label{prop:G_rho}
The penalized lifted regularizer $G_{\boldsymbol{\rho}}$ satisfies the
following properties:
\begin{enumerate}
    \item \textbf{Lower approximation.} $G_{\boldsymbol{\rho}}(x)
    \leq R(x)$ for all $x$ and all $\boldsymbol{\rho} \geq 0$.

    \item \textbf{Monotonicity and concavity.} If
    $\boldsymbol{\rho}_1 \leq \boldsymbol{\rho}_2$ componentwise,
    then $G_{\boldsymbol{\rho}_1}(x) \leq G_{\boldsymbol{\rho}_2}(x)$
    for all $x$. In the scalar-parameter setting
    $\rho_1 = \cdots = \rho_L = \rho$, the map $\rho \mapsto
    G_\rho(x)$ is nondecreasing and concave for every fixed $x$.

    \item \textbf{Finite exactness.} There exists a finite threshold
    $\bar{\boldsymbol{\rho}}$ such that $G_{\boldsymbol{\rho}}(x) =
    R(x)$ for all $x$ whenever $\boldsymbol{\rho} \geq
    \bar{\boldsymbol{\rho}}$. In the scalar-parameter setting, an
    explicit threshold is
    \begin{equation}\label{eq:rho_bar}
        \bar{\rho} = \max_{\ell=1,\ldots,L}
        \Bigl\{
            \|V\|_\infty \, L_\sigma^{L-\ell+1}
            \prod_{j=\ell+1}^{L+1} \|W_j\|_{1\to1}
        \Bigr\},
    \end{equation}
    where $L_\sigma$ is the Lipschitz constant of $\sigma$.

    \item \textbf{Exactness for monotone convex networks.} If
    $\sigma$ is nondecreasing and $W_\ell \geq 0$ componentwise
    for all $\ell$, then $G_{\boldsymbol{\rho}}(x) = R(x)$ for
    all $x$ and all $\boldsymbol{\rho} \geq 0$.
\end{enumerate}
\end{proposition}

Property~1 shows that $G_{\boldsymbol{\rho}}$ is always a lower approximation of $R$, tightening monotonically as $\rho$ increases (Property~2). Property~3 guarantees that the original problem is recovered exactly at a finite and computable penalty. Property~4 recovers the known exactness of the unpenalized lifting for ICNNs as a special case. We note that this result concerns a particular network architecture rather than an arbitrary convex function $R$: the construction of $G_{\boldsymbol{\rho}}$ depends on the lifted representation induced by the network, not only on $R$ itself. Consequently, different network architectures representing the same function $R$ may lead to different relaxations $G_{\boldsymbol{\rho}}$ with different exactness thresholds $\bar\rho$. Together, these properties define a natural homotopy
\begin{equation}\label{eq:homo-Grho}
    D+\lambda G_{\rho_1} \leq D+\lambda G_{\rho_2} \leq \cdots \leq D+\lambda G_{\rho_K} = D+\lambda R
\end{equation}
along an increasing schedule $\rho_1 < \rho_2 < \cdots < \rho_K =
\bar\rho$, from a loose lower approximation at small $\rho$ to the
exact inference problem at $\rho = \bar\rho$, which we exploit
algorithmically in Section~\ref{sec:optm}. The concavity of $\rho\mapsto G_\rho(x)$ (Property~2) implies diminishing returns as $\rho$ increases, taht is the approximation
tightens quickly at first and more slowly near $\bar\rho$, suggesting a geometrically spaced $\rho$ schedule. Proofs are given in the supplementary material.

\subsection{DC decomposition}\label{sec:dc-decomp}
The continuation path defined by $\{G_{\boldsymbol{\rho}}\}$
is useful algorithmically only if each lifted problem can be
solved efficiently. The following proposition shows that, under
convex activations, each problem admits a DC decomposition.

\begin{proposition}
\label{prop:DC_convex_activation_L}
Assume that $\sigma:\mathbb{R}\to\mathbb{R}$ is convex and
applied componentwise, and that $D$ is convex. Then the
optimization problem
\begin{equation}
\label{eq:min_g_G_L}
\min_{x} \; \left\{D(x) + \lambda G_{\boldsymbol{\rho}}(x) \right\}
\end{equation}
is equivalent to the lifted DC problem
\begin{equation}
\label{eq:lifted_DC_problem_L}
\min_{x,\boldsymbol{z}} \;
P_{\lambda,\boldsymbol{\rho}}(x,\boldsymbol{z})
-
Q_{\lambda,\boldsymbol{\rho}}(x,\boldsymbol{z}),
\end{equation}
where
\begin{align}
P_{\lambda,\boldsymbol{\rho}}(x, \boldsymbol{z})
&=
D(x)
+
\lambda V^\top z_L
+
\lambda
\sum_{\ell=1}^{L}
\rho_\ell \mathbf{1}^\top z_\ell
+
\iota_{\mathcal{C}_{\sigma}^{L}},
\label{eq:P_DC_L}
\\
Q_{\lambda,\boldsymbol{\rho}}(x, \boldsymbol{z})
&=
\lambda \sum_{\ell=1}^{L}
\rho_\ell \mathbf{1}^\top
\sigma(\phi_\ell).
\label{eq:Q_DC_L}
\end{align}
Both $P_{\lambda,\boldsymbol{\rho}}$ and
$Q_{\lambda,\boldsymbol{\rho}}$ are convex, and
$\mathcal{C}_{\sigma}^{L}$ is a convex set.
\end{proposition}

In the case where $\sigma$ is piecewise linear, the DC formulation allows the problem to be decomposed into a sequence of quadratic programs.

\begin{corollary}[DC formulation with linear constraints]
\label{cor:DC_piecewise_affine}
If $\sigma$ is additionally piecewise affine, then
$\mathcal{C}_{\sigma}^{L}$ is polyhedral and each DCA
subproblem reduces to a quadratic program (QP) with linear
constraints.
\end{corollary}

The convexity of $P_{\lambda,\boldsymbol{\rho}}$ and $Q_{\lambda,\boldsymbol{\rho}}$ follows respectively from the convexity of $D$, the linearity of the lifted terms, and the convexity of $\iota_{\mathcal{C}_\sigma^L}$, and from the fact that each term $\rho_\ell\mathbf{1}^\top \sigma(\phi_\ell)$ is a composition of a convex function with an affine map. The key consequence is that each DCA step requires solving a \emph{globally} convex subproblem. For piecewise-affine activations such as ReLU (Corollary~\ref{cor:DC_piecewise_affine}), this subproblem is a QP. For general convex activations such as smoothed ReLU, it can be formulated as a saddle-point problem addressable by PDHG \cite{chambolle2011first}.

\begin{remark}[Practical evaluation of $G_\rho$]\label{rem:practicaleval}

For a fixed $x$, the piecewise-linear structure of $\sigma$ allows
$G_{\rho}(x)$ to be expressed as the minimum of a collection of linear
programs, one for each possible activation pattern, corresponding to
the different linear pieces of $\sigma$. If $
m=\sum_{\ell=2}^L n_\ell
$ denotes the number of binary variables defining these patterns, this amounts
to solving up to $2^m$ linear programs for each $x$.

In practice, we obtained good approximation results exploiting the activation pattern $\varepsilon_R(x)$ obtained
from the forward pass of the original network. The optimal pattern for $G_\rho(x)$
is typically found within a small Hamming neighborhood of
$\varepsilon_R(x)$. Restricting the search to patterns satisfying
$d_H(\varepsilon,\varepsilon_R(x))\le r$ reduces the number of LPs to
$
\sum_{k=0}^r \binom{m}{k},
$
which is polynomial in $m$ for fixed $r$, instead of exponential. 
This strategy is currently heuristic, as no theoretical guarantee on the
distance between $\varepsilon_R(x)$ and an optimal pattern is available.
\end{remark}

\section{Homotopic-DC Algorithm}
\label{sec:optm}

\subsection{Homotopic continuation}\label{sec:homo-cont}

The continuation path established in Section~\ref{sec:lifting} naturally suggests a homotopic optimization strategy. Denoting the full lifted variable $v = (x, z_1,
\ldots, z_L)$, at each level $k$ we approximately solve
\[
    v^{(k)} \approx \arg\min_{v} \bigl\{
        P_{\lambda,\rho_k}(v) - Q_{\lambda,\rho_k}(v)
    \bigr\}
\]
initialized at $v^{(k-1)}$.

The resulting sequence
\[
    v^{(0)} \longrightarrow v^{(1)}
    \longrightarrow \cdots
    \longrightarrow v^{(K)}
\]
tracks a solution path as the relaxation is progressively tightened. Warm-starting improves the effectiveness of the scheme: at each level $k$, the solution $v^{(k-1)}$ from the previous level is already close to feasible for the tighter constraints at $\rho_k$, providing a better initialization than a cold start and allowing DCA to converge in fewer iterations.

Since $G_{\bar\rho} = R$, the final stage corresponds exactly to the original inference problem \ref{eq:P}, so the homotopy terminates at a solution of the original nonconvex problem rather than an approximation. At small $\rho$, the lifted problem is a loose relaxation with a simpler landscape, allowing the iterates to find a good basin of attraction. As $\rho$ increases, the constraints tighten progressively, steering the iterates toward a high-quality solution of the original problem.

\begin{figure}[t]
\centering
\includegraphics[width = 0.9\linewidth]{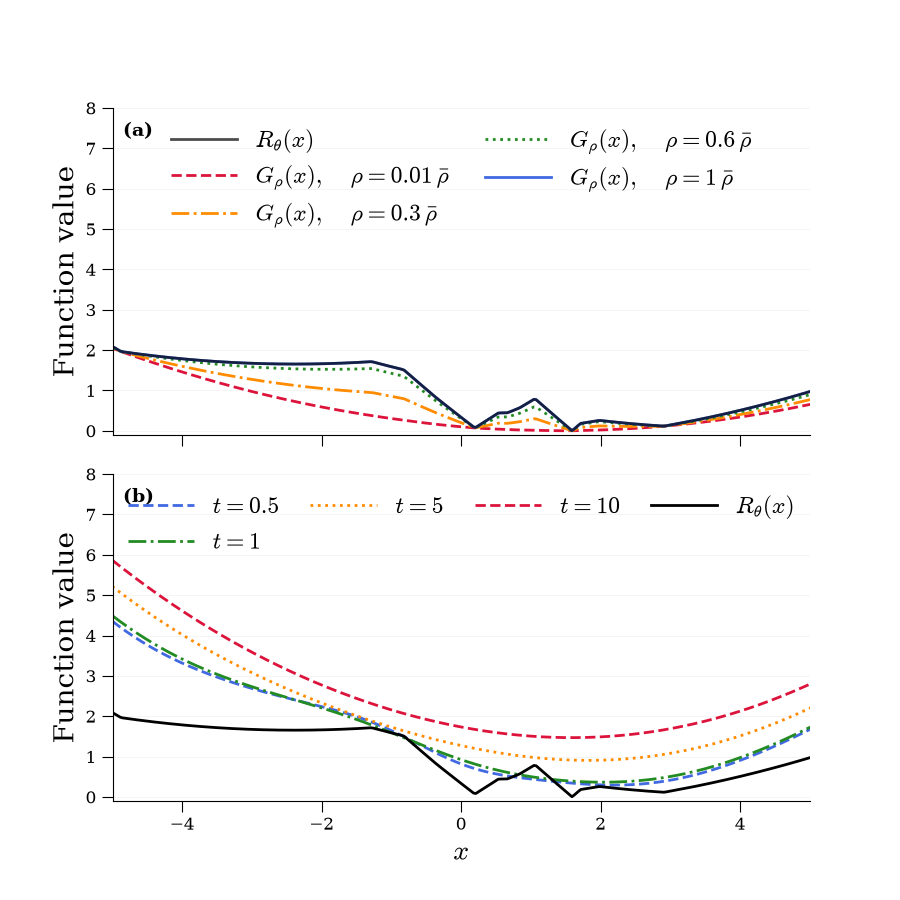}
\caption{Visual comparison between the homotopy defined in~\ref{eq:homo-Grho} and the Gaussian homotopy, for a function of the form the objective value of ~\eqref{eq:P}, where $R$ is a neural network 
randomly generated as described in Section~\ref{sec:lowdim}. We observe the properties 
established in Proposition~\ref{prop:G_rho}: $G_\rho$ always lies below the target 
function, and $G_\rho$ is monotonically increasing with respect to $\rho$. }\label{fig:gaussDC}
\end{figure}

This strategy is an instance of \emph{continuation}, or \emph{homotopy optimization}: instead of attacking the target nonconvex objective directly, one constructs a parameterized family of problems, starts from a simpler optimization landscape, and progressively deforms it toward the original objective while warm-starting each problem from the solution of the preceding one~\cite{
iwakiri2022single,sato2025explicit}. A prominent example is \emph{Gaussian homotopy}, where the objective is progressively smoothed through convolution with Gaussian kernels. Such graduated optimization strategies have been shown empirically to reach better stationary points than direct local optimization in a variety of nonconvex problems. Our approach differs in that the continuation path is \emph{network-intrinsic}, in the sense that it is induced by the internal architecture rather than external smoothing. Moreover, it lies \emph{below} $R$ at all $\rho$ (Figure~\ref{fig:gaussDC}a), preserving the local structure of the landscape, in contrast to Gaussian homotopy which may lie above $R$ and smooths away local structure (Figure~\ref{fig:gaussDC}b). Furthermore, it comes with the additional guarantee of exact recovery at a finite and computable $\bar\rho$.

\subsection{DCA for each subproblem}\label{sec:dca}

Each lifted DC problem along the path is solved using the
\emph{Difference-of-Convex Algorithm} (DCA)~\cite{lethi2012exact}.
Given the current iterate $v^k$, DCA linearises
$Q_{\lambda,\rho}$ at $v^k$ and solves the resulting convex
subproblem:
\begin{equation}
    v^{k+1} = \arg\min_{v} \bigl\{
        P_{\lambda,\rho}(v)
        - \langle \nabla Q_{\lambda,\rho}(v^k),\, v \rangle
    \bigr\}.
    \label{eq:DCA_sub}
\end{equation}
DCA guarantees monotone descent on the lifted objective:
\[
    P_{\lambda,\rho}(v^{k+1}) - Q_{\lambda,\rho}(v^{k+1})
    \leq
    P_{\lambda,\rho}(v^{k}) - Q_{\lambda,\rho}(v^{k}).
\]
For piecewise-affine activations, the subproblem~\eqref{eq:DCA_sub} reduces to a QP with linear
constraints (Corollary~\ref{cor:DC_piecewise_affine}), which can be solved via an off-the-shelf QP solver for small-scale problems. For large-scale problems, the subproblem can be formulated as a saddle-point problem, which can be addressed efficiently by
PDHG~\cite{chambolle2011first}.
The full procedure is summarized in Algorithm~\ref{alg:hdc}.

\begin{algorithm}[t]
\caption{H-DC}\label{alg:hdc}
\begin{algorithmic}[1]
\Require Initial point $v^{(0)}$, schedule
         $\rho_1 < \cdots < \rho_K = \bar\rho$,
         max DCA iterations $T$
\For{$k = 1, \ldots, K$}
  \State $v \leftarrow v^{(k-1)}$
  \For{$t = 1, \ldots, T$}
    \State Compute $\nabla Q_{\lambda,\rho_k}(v)$
    \State $v_{\mathrm{prev}} \leftarrow v$
    \State $v \leftarrow \arg\min_{v'} \bigl\{
           P_{\lambda,\rho_k}(v')
           - \langle \nabla Q_{\lambda,\rho_k}(v),\, v' \rangle
           \bigr\}$
    \If{$\|v - v_{\mathrm{prev}}\| \leq \texttt{tol}$}
      \textbf{break}
    \EndIf
  \EndFor
  \State $v^{(k)} \leftarrow v$
\EndFor
\Return $v^{(K)}$
\end{algorithmic}
\end{algorithm}

\section{Numerical experiments}
\label{sec:numerics}

We evaluate H-DC on three settings of different complexity. First, we consider \emph{randomly generated
inverse problems} in one and two dimensions with randomly generated nonconvex regularizers, where the
global minimum can be computed to a high accuracy and we can measure the
optimality gap directly. Second, we consider a \emph{manifold denoiser regularizer} on the unit circle
$S^1$, where the signal class has known geometry and we can visualize the optimization trajectories. Third, we consider \emph{image inpainting} on CIFAR-10 with
a learned approximation of the Mumford-Shah regularizer, where the quality of the local minimum directly affects reconstruction quality. In each setting, H-DC is
compared against gradient-based method initialized from the same
starting point. All experiments use piecewise-affine activations. DCA subproblems are solved via QP for the first two experiments and via
PDHG for the image experiment.

\subsection{Inverse problems in low dimension}\label{sec:lowdim}
We consider randomly generated inverse problems in dimensions one and two. The data-fidelity term is
\[
    D(x) = \|Ax - b\|^2,
\]
where $A$ and $b$ are sampled from $\mathcal{N}(0,1)$ and $A$ is subsequently normalized. The regularizer $R$ is represented by a neural network with absolute-value
activation $\sigma(t) = |t|$, which is convex and Lipschitz, so Proposition~\ref{prop:G_rho} applies and the explicit threshold $\bar\rho$ is available.

\paragraph{Dimension 1.}
We use two hidden layers of width $5$. The entries of $W_1$ are sampled from $\mathcal{N}(0,1)$, while the entries of $W_\ell$ and $W_{\mathrm{out}}$ are sampled
from Gaussian distributions with standard deviation $1/\sqrt{n_{\ell-1}}$ and $1/\sqrt{n_L}$ respectively, where $n_\ell$ denotes the width of layer $\ell$. The biases are sampled from $\mathcal{U}([-1,1])$ and $V$ from $\mathcal{U}([0,1])$. To ensure coercivity, we retain only networks whose asymptotic slopes in both
directions are at least $0.1$, since $R$ is asymptotically linear as $x \to \pm\infty$.

\paragraph{Dimension 2.}
We consider $R : \mathbb{R}^2 \to \mathbb{R}$ with two hidden layers of width $3$ to $5$ and activation
$\sigma(t) = |t|$. The initialization follows the same
rules as the one-dimensional setting, with weight scaling adapted to the input dimension and $A_\ell\sim \mathcal{N}(0, 1/2)$.

%The entries of each weight matrix are independently sampled from
%$(W_\ell)_{ij}\sim
%\mathcal{N}\left(0,\frac{1}{n_{\ell-1}}\right),$
%where $n_{\ell-1}$ denotes the dimension of the input to layer $\ell$.
% The matrices $A_\ell$ are sampled from
%$\mathcal{N}(0,1/d)$. The biases are independently sampled
%from $\mathcal{U}(-1,1)$ and the entries of $V$ from
%$\mathcal{U}(0.5,1)$. The generated networks are subsequently screened,
%and only instances satisfying the prescribed nonconvexity criteria are
%retained for the numerical experiments.

\begin{figure}[t]
\centering
\includegraphics[width=\linewidth]{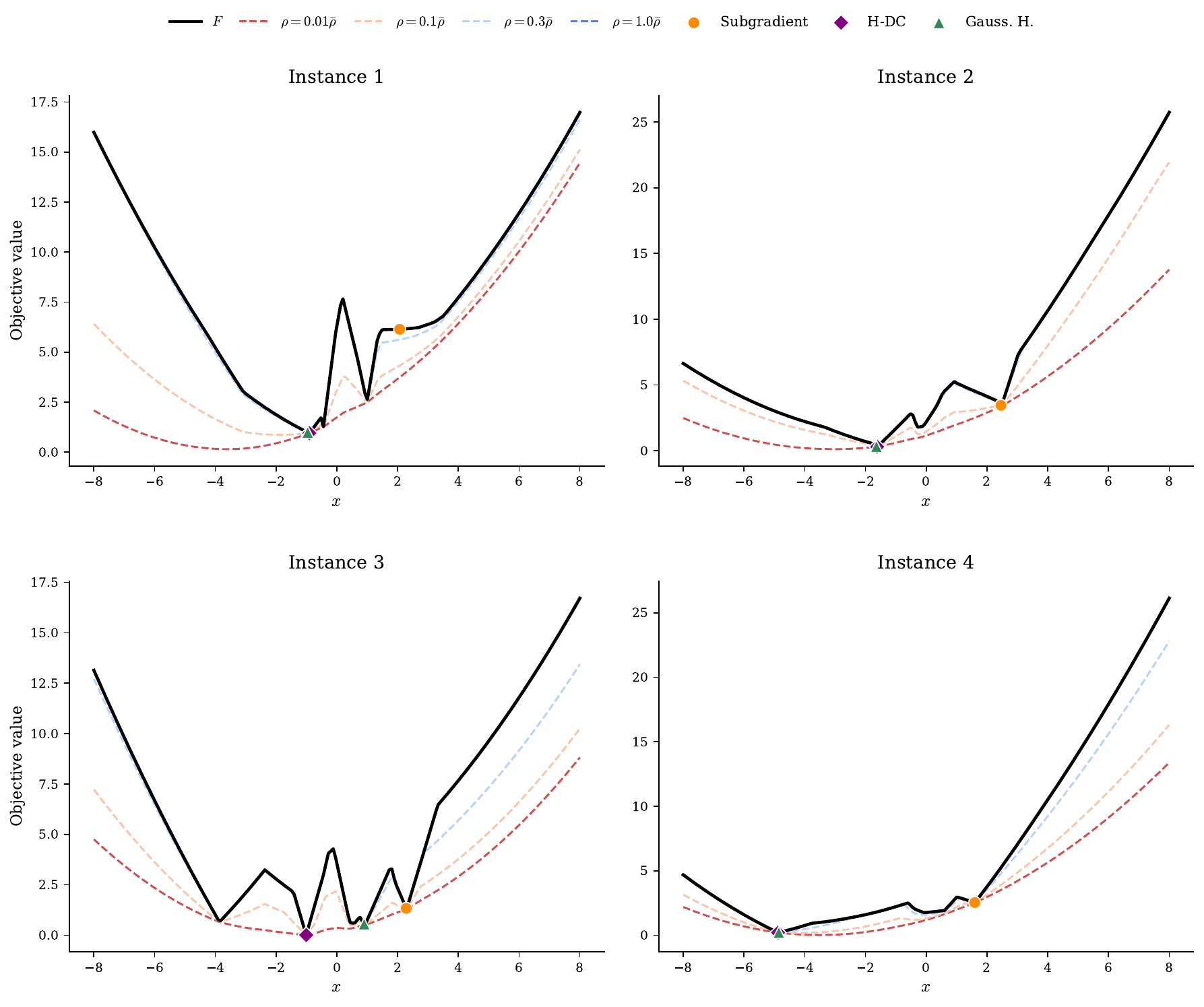}
\caption{Four randomly generated 1D instances. The penalized lifted objectives $G_\rho$ (dashed) and oribinal objective (black) are shown, along with solutions of H-DC (purple diamond), subgradient (orange circle) and Gaussian homotopy (green triangle).}
\label{fig:instances}
\end{figure}

\paragraph{Results.}
Figure~\ref{fig:instances} shows four 1D instances. H-DC outperforms subgradient on all four, reaching solutions with lower objective values. Gaussian homotopy performs competitively, matching H-DC on three instances but converging to a suboptimal point on instance 3. The homotopy path $\{G_\rho\}$ lies below $F$ at all $\rho$ and converges to it as $\rho \to \bar\rho$, directly confirming the theoretical properties of Proposition~\ref{prop:G_rho}. 

\begin{figure}
\centering
\includegraphics[width = 0.9\linewidth]{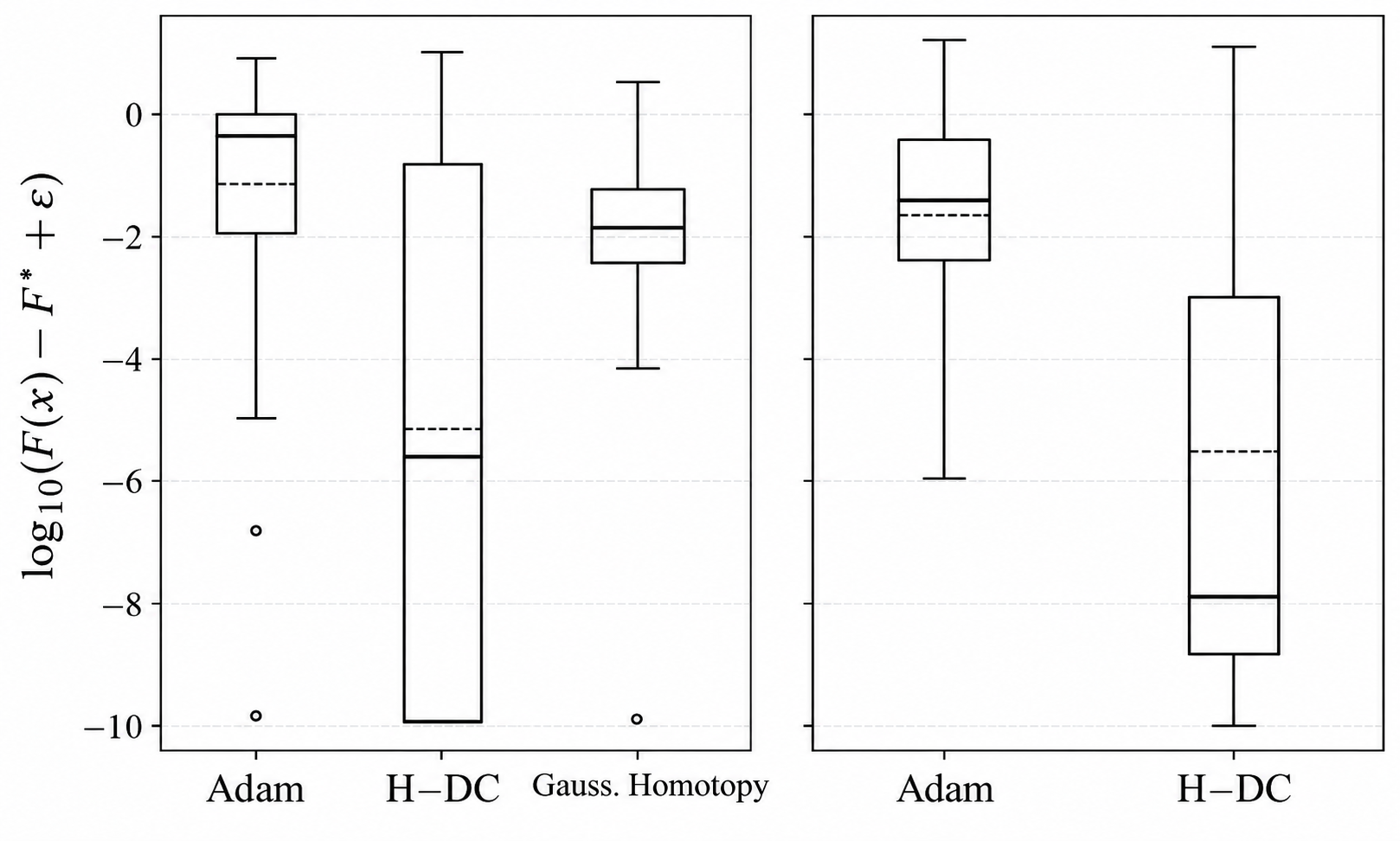}
\caption{\textbf{Right:} Comparison of the optimality gap obtained with Adam and H-DC and Gaussian homotopy in dimension one. \textbf{Left:} Same comparison in dimension two. The median is shown by a solid black line, while the mean is shown by a dashed line. }\label{fig:boxplots}
\end{figure}

Figure~\ref{fig:boxplots} extends this comparison to a large number of randomly
generated problems. The mean optimality gap decreases from approximately $0.2$ for Adam to $10^{-5.6}$ for H-DC, a reduction of several orders of magnitude, and H-DC also outperforms Gaussian homotopy on average across both dimensions. These results confirm that by exploiting the epigraph structure of the network, H-DC consistently escapes poor local minima that trap gradient-based and smoothing-based methods.

\subsection{Manifold denoiser regularizer}\label{sec:mfd-denoise}
\begin{figure*}[t]
\centering
\includegraphics[width=\linewidth]{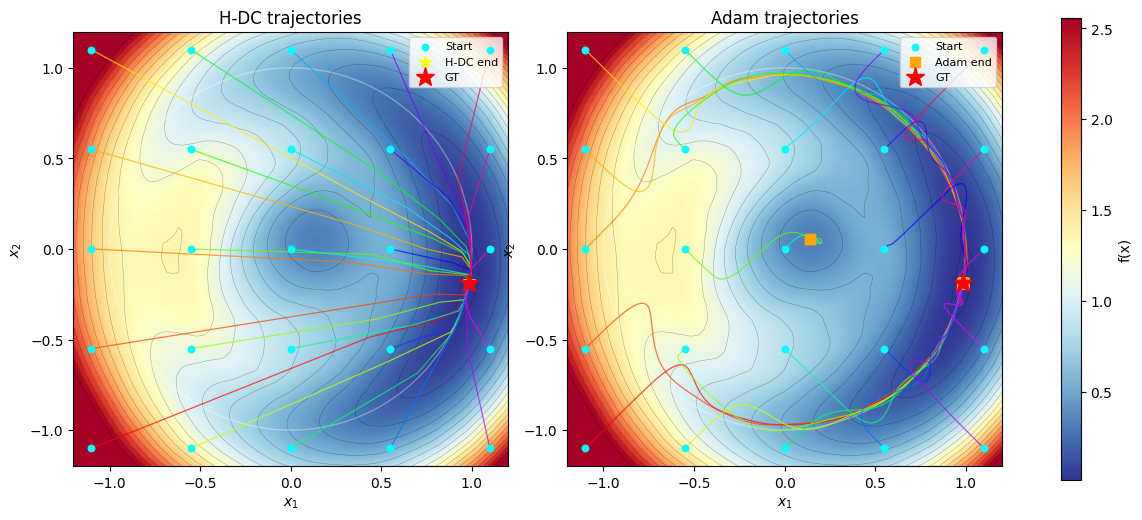}
\caption{Optimization trajectories from a $5\times 5$
grid of initializations for H-DC (left) and Adam (right). Ground truth $x^*$ marked by a red star.}
\label{fig:circle}
\end{figure*}

We consider a two-dimensional Gaussian sensing problem
where the signal is assumed to lie on the unit circle
$S^1 \subset \mathbb{R}^2$. The data fidelity is $D(x) = \tfrac{1}{2}\|Ax - y\|^2$ with $A \in\mathbb{R}^{3 \times 2}$ random normalized and measurements $y = Ax^* + \varepsilon$, where $x^* \in
S^1$ and $\varepsilon \sim \mathcal{N}(0, \sigma^2 I)$,
$\sigma = 0.1$. A three-layer dense denoiser $N_\theta:\mathbb{R}^2 \to \mathbb{R}^2$ with ReLU activations
and hidden width $128$ is trained to map noisy observations back to $S^1$ by minimizing $\mathbb{E}[\|N_\theta(x + \varepsilon) - x\|^2]$ with $x \sim \mathrm{Uniform}(S^1)$ and $\varepsilon \sim
\mathcal{N}(0, \sigma^2 I)$, $\sigma = 0.1$, for 5000
epochs using AdamW with learning rate $10^{-3}$ and batch size $256$. We consider the LSR regularizer
\[
    R(x) = \tfrac{\alpha}{2}\|x - N_\theta(x)\|^2,
\]
which penalizes $x$ for lying far from the fixed-point set of $N_\theta$, an approximation of $S^1$. Although the LSR regularizer can be cast within architecture~\eqref{eq:network} as described in Remark~1, in practice we denote $z_L = N_\theta(x) =
W_L z_{L-1} + b_L$ the linear output of the denoiser, and rewrite the full objective as
\begin{align*}
    \min_{x, z_1, \ldots, z_{L-1}} \;
    &\tfrac{1}{2}\|Ax - y\|^2
    + \tfrac{\alpha}{2}\|x - W_L z_{L-1} - b_L\|^2 \\
    &+ \lambda \rho \sum_{\ell=1}^{L-1}
    \mathbf{1}^\top (z_\ell - \sigma(\phi_\ell)) \\
    \text{s.t.} \quad
    &z_\ell \geq \sigma(\phi_\ell),
    \quad \ell = 1, \ldots, L-1,
\end{align*}
where the quadratic residual is treated as part of the convex component $P_{\lambda,\boldsymbol{\rho}}$, and lifting and penalization are applied only to the hidden ReLU layers $z_1, \ldots, z_{L-1}$. The subproblem remains a QP with linear constraints, and is solved via the off-the-shelf QP solver cvxopt~\cite{andersen2013cvxopt}.

We run H-DC and Adam from a $5 \times 5$ grid of initializations covering $[-1.5, 1.5]^2$ with regularization parameter $\alpha=10$, and record the final objective value and trajectory of each method. The $\rho$ schedule is set to $[0.001, 0.005, 0.01, 0.05, 0.1]$.

\begin{remark}
The LSR regularizer $R(x) = \frac{\alpha}{2}\|x -N_\theta(x)\|^2$ does not satisfy the Lipschitz assumption of Proposition~\ref{prop:G_rho}, as the quadratic activation $\sigma(s) = s^2$ is not globally Lipschitz. In particular, the finite exactness result and the explicit threshold $\bar\rho$ in~\eqref{eq:rho_bar} do not necessarily apply. However, as shown above, the quadratic term is handled directly as part of $P_{\lambda,\boldsymbol{\rho}}$ and is never penalized, so the subproblem remains a well-posed QP. H-DC is applied with an manually chosen $\rho$ schedule and is shown empirically to reach the global minimum from all initializations.
\end{remark}

\paragraph{Results.}
As shown in Figure~\ref{fig:circle}, H-DC converges to the global minimum from all 25 initializations, including those near the origin where Adam stagnates. The trajectory plot reveals a qualitative difference between the two methods: H-DC navigates directly toward the ground truth as $\rho$ increases along the homotopy path, while Adam follows the local gradient and terminates at a spurious local minimum. This illustrates the key advantage of the homotopic scheme: by starting from a loose relaxation and progressively tightening, H-DC escapes basins of attraction that trap gradient-based methods.

\subsection{Image inpainting with a learned regularizer}

Natural images exhibit rich local structure such as sharp edges and smooth regions. The truncated quadratic regularizer~\cite{strekalovskiy2014real,deng2021preconditioned}
\begin{equation}
    R_{\mathrm{MS}}(x) = \sum_{i}
    \min\!\left\{\alpha \|\nabla_i x\|^2,\, \beta\right\},
    \label{eq:ms}
\end{equation}
a discretization of the Mumford-Shah functional~\cite{mumford1989optimal}, promotes piecewise smooth reconstructions by saturating the penalty for large gradients, where $\nabla_i x$ denotes the image gradient at pixel $i$, $\alpha = 20$ and $\beta = 0.1$. This regularizer is nonconvex and known to be NP-hard to minimize in general~\cite{alexeev2010complexity}, making it a challenging but practically relevant test case.

We approximate $R_{\mathrm{MS}}$ using a fully connected network $R$ following architecture~\eqref{eq:network} with ReLU activations. Since $R_{\mathrm{MS}}$ is a sum of a scalar potential $\psi(t) = \min\{\alpha t^2,\beta\}$ applied to image gradient magnitudes, we train the network to approximate the scalar potential $\psi$ directly, and apply it to the gradient magnitude $\|\nabla_i x\|$ at each pixel. The network consists of three hidden layers of width $24$ with skip connection. The network is trained on CIFAR-10 training
images by minimizing the mean squared error between $N(t)$ and $\psi(t)$ using AdamW with learning rate $10^{-3}$ for $100$ epochs and batch size $256$. The learned regularizer is then $R(x) = \sum_iN(\|\nabla_i x\|)$, which inherits the piecewise-affine structure of the ReLU network. The DCA subproblem is solved via PDHG, with a $\rho$ schedule $[0.02, 0.1, 0.3, 1.0] \cdot \bar\rho$, where $\bar\rho$ is computed via~\eqref{eq:rho_bar}. 

\begin{figure*}[t]
\centering
\includegraphics[width=\linewidth]{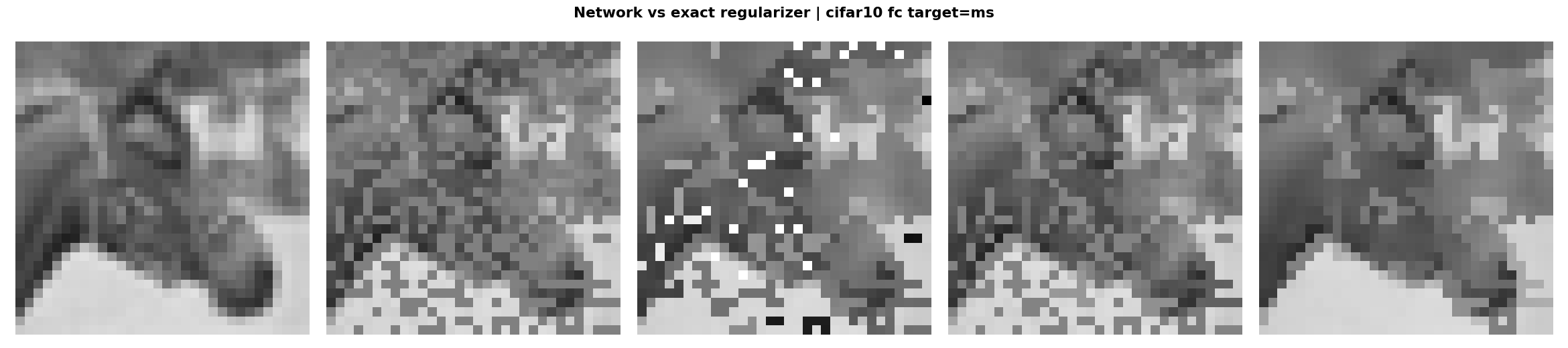}
\caption{Inpainting results (30\% pixels kept). From
left to right: ground truth, observation, Adam with
learned $R$ (PSNR=15.4\,dB), Adam with exact
$R_{\mathrm{MS}}$ (PSNR=19.2\,dB), H-DC with learned
$R$ (PSNR=23.6\,dB).}
\label{fig:inpainting}
\end{figure*}

We consider inpainting on CIFAR-10 test images, where $70\%$ of pixels are randomly kept and the remainder
masked. Reconstruction is obtained by solving
\[
    \min_x \; \tfrac{1}{2}\|M \odot (x - y)\|^2
    + \lambda R_\theta(x),
\]
where $M$ is the binary mask and $y$ is the observed image, with $\lambda = 0.02$. We compare H-DC against Adam initialized from the same starting point, running Adam until full convergence to ensure a fair comparison. To contextualize the role of the network approximation, we also compare against Adam with the exact regularizer~\eqref{eq:ms} on a representative example.

\begin{figure}[t]
\centering
\includegraphics[width=\linewidth]{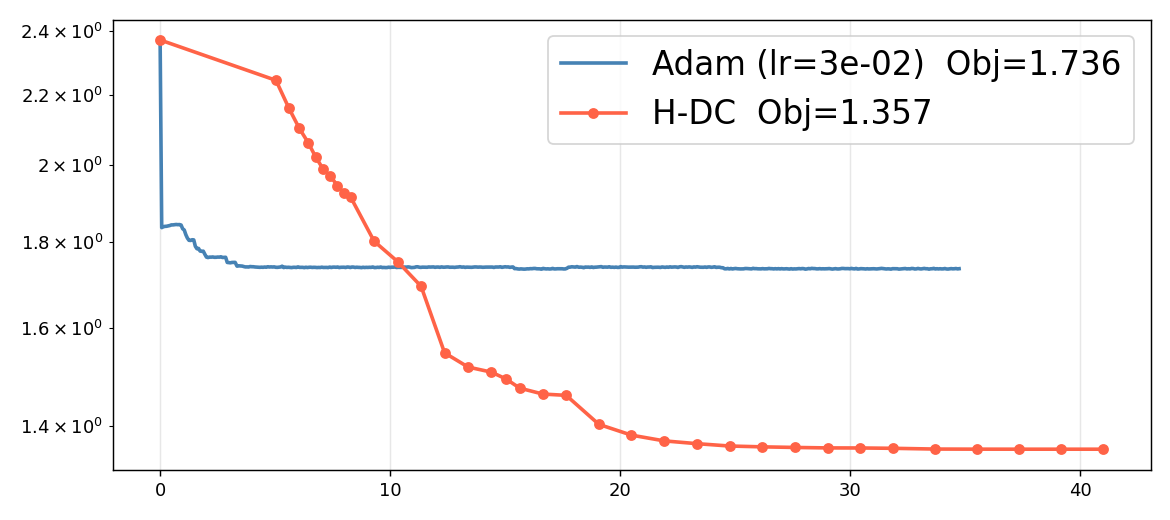}
\caption{Objective value vs.\ time for H-DC and Adam
with learned $R_\theta$ ($\lambda=0.02$).}
\label{fig:ms_convergence}
\end{figure}

\paragraph{Results.}
Figure~\ref{fig:inpainting} shows a sample reconstruction. Adam with the learned regularizer leaves visible inpainting artifacts (PSNR=15.4\,dB), while H-DC recovers finer image structure (PSNR=23.6\,dB). Adam with the exact $R_{\mathrm{MS}}$ achieves 19.2\,dB, an improvement over Adam with the learned regularizer, but still below H-DC, suggesting that the gain comes from better optimization rather than regularizer quality.

Figure~\ref{fig:ms_convergence} shows the objective value vs.\ time for H-DC and Adam with the learned regularizer. Since each H-DC step involves solving a convex subproblem while Adam takes a single gradient step, iteration counts are not directly comparable. We therefore report wall-clock time and run Adam until full convergence. H-DC continues to decrease the objective well below the level where Adam stagnates.

Table~\ref{tab:inpainting} confirms these findings over 10 CIFAR-10 test images. H-DC reaches a lower objective value than Adam on all 10 images, with a mean gap of $0.28$ and a mean PSNR improvement of $6.5$\,dB ($23.53 \pm 1.81$ vs.\ $17.00 \pm 1.72$\,dB).

\begin{table}[t]
\centering
\caption{Average results over 10 CIFAR-10 test images,
inpainting with 30\% pixels kept, learned $R$
($\lambda=0.02$).}
\label{tab:inpainting}
\begin{tabular}{lcc}
\toprule
Method & Obj $\downarrow$ & PSNR (dB) $\uparrow$ \\
\midrule
Adam & $1.47 \pm 0.25$ & $17.00 \pm 1.72$ \\
H-DC & $\mathbf{1.19 \pm 0.24}$ &
    $\mathbf{23.53 \pm 1.81}$ \\
\bottomrule
\end{tabular}
\end{table}

\section{Conclusion}
We introduced H-DC, a homotopic optimization scheme for network inference problems. By lifting the network into an epigraph formulation and progressively tightening the relaxation via an increasing penalty parameter $\rho$, H-DC solves a sequence of DC subproblems that provably recover the original objective at a finite and computable threshold $\bar\rho$. Each subproblem is solved globally via DCA, reducing to a QP for piecewise-affine activations or a saddle-point problem addressed by PDHG.

Experiments on synthetic problems, manifold denoiser regularizer, and image inpainting with learned regularizer demonstrate that H-DC consistently reaches solutions with lower objective values than gradient-based methods. Notably, H-DC with a learned approximation outperforms Adam with the exact regularizer, highlighting that optimization quality matters as much as regularizer quality.

Several directions remain open. Extending the finite exactness result to more general activations, scaling to larger networks, and developing adaptive $\rho$ schedules with convergence guarantees are natural
directions.

\newpage
% \subsection*{AI use statement}

% (This section is \textbf{required} and does not count toward the page
% limit.)

% % TODO (author): complete this disclosure honestly before submission.
% \subsection*{AI use statement}
% In this work, we used generative AI tools to aid in implementing the methodology, executing comparison methods. We did not use generative AI tools to develop theoretical models or conceptual frameworks, to formulate mathematical claims, provide critical ingredients for proving mathematical claims, to propose or refine hypotheses, or to interpret the results. Synthetic datasets, assisting with translation, cleaning and reformatting datasets, and supporting qualitative and thematic data analysis are not applicable to this work. Additionally, we used generative AI tools for visualizing the results of our method, improving the written clarity and readability of the manuscript including tables and algorithms. We have reviewed all AI-assisted work, including written text, code, visualizations, and technical details. We take responsibility for the final content of this work, including text, claims, or artifacts produced with the aid of generative AI.

% See the AISTATS 2027 AI Policy for Authors for more details. This
% statement should not be more than 1 page.

%\subsubsection*{Acknowledgements}
% TODO (author): acknowledgements go here in the camera-ready version
% only, to preserve anonymity during review.

\bibliographystyle{unsrtnat}
\bibliography{iclr2027_conference}

\clearpage
\appendix
\thispagestyle{empty}

% Supplementary material: single-column format, per AISTATS instructions.
\onecolumn
\aistatstitle{Taming Nonconvexity with Epigraph Lifting: \\
Supplementary Materials}

\section{Proof of Proposition 1}

\paragraph{1. Lower approximation.}
Let $x\in X$ and consider the variables obtained by the forward pass,
\[
z_1=\sigma(\phi_1(x)),~~
z_\ell=\sigma(\phi_\ell(x,z_{\ell-1})),\quad \ell=2,\ldots,L,
\]
and
\[
u=\sigma(W_{L+1}z_L+b_{L+1}).
\]
This defines a feasible point of $\mathcal C_\sigma^L$ for which all
penalty terms vanish. Hence, by the definition of $G_{\boldsymbol\rho}$,
\[
G_{\boldsymbol\rho}(x)\leq V^\top u=R(x).
\]

\paragraph{2. Monotonicity and concavity with respect to the penalty parameter.}
Let $\boldsymbol{\rho}_1\leq \boldsymbol{\rho}_2$ componentwise. Since
\[
z_\ell-\sigma(\phi_\ell(x,z_{\ell-1}))\geq 0,
\qquad \ell=1,\ldots,L,
\]
for every $(x,z_1,\ldots,z_L,u)\in\mathcal C_\sigma^L$, we have

\[
V^\top u
+\sum_{\ell=1}^L \rho_{1,\ell}\,
\mathbf 1^\top
\bigl(z_\ell-\sigma(\phi_\ell(x,z_{\ell-1}))\bigr)
\leq
V^\top u
+\sum_{\ell=1}^L \rho_{2,\ell}\,
\mathbf 1^\top
\bigl(z_\ell-\sigma(\phi_\ell(x,z_{\ell-1}))\bigr).
\]
Taking the infimum over the same feasible set yields
\[
G_{\boldsymbol{\rho}_1}(x)
\leq
G_{\boldsymbol{\rho}_2}(x).
\]

In the scalar-parameter setting, for fixed $x$ and fixed feasible
$(z_1,\ldots,z_L,u)$, the objective defining $G_\rho(x)$ is affine in
$\rho$. Hence, $G_\rho(x)$ is the pointwise infimum of affine functions
of $\rho$, and is therefore concave. Its monotonicity follows from the
first part.

\paragraph{3 and 4. Explicit exactness threshold.}
Assume that $\sigma$ is $L_\sigma$-Lipschitz continuous. 
We show that the relaxation constraints are active at any minimizer 
whenever the penalty parameters are sufficiently large.

Consider first the last hidden layer. A perturbation 
$\delta z_L$ induces a variation of the output satisfying

\[
    \left|
V^\top \sigma(W_{L+1}(z_L+\delta z_L)+b_{L+1})
-
V^\top \sigma(W_{L+1}z_L+b_{L+1})
\right|
\leq
\|V\|_\infty
L_\sigma
\|W_{L+1}\delta z_L\|_1
\\
\leq
\|V\|_\infty
L_\sigma
\|W_{L+1}\|_{1\to1}
\|\delta z_L\|_1.
\]
Hence, if
\[
\rho_L
\geq
\|V\|_\infty
L_\sigma
\|W_{L+1}\|_{1\to1},
\]

any positive relaxation gap at layer $L$ can be reduced without
increasing the objective. Thus, an optimal solution can be chosen such
that
\[
z_L=\sigma(\phi_L(x,z_{L-1})).
\]

Proceeding backward, a perturbation $\delta z_\ell$ at layer $\ell$
propagates through the subsequent layers according to
\[
\|\delta z_j\|_1
\leq
L_\sigma
\|W_j\|_{1\to1}
\|\delta z_{j-1}\|_1,
\qquad j=\ell+1,\ldots,L+1.
\]
Consequently, its effect on the output is bounded by
\[
\|V\|_\infty
L_\sigma^{L-\ell+1}
\left(
\prod_{j=\ell+1}^{L+1}
\|W_j\|_{1\to1}
\right)
\|\delta z_\ell\|_1.
\]
It is therefore sufficient to choose
\[
\bar\rho_\ell
=
\|V\|_\infty
L_\sigma^{L-\ell+1}
\prod_{j=\ell+1}^{L+1}
\|W_j\|_{1\to1},
\qquad
\ell=1,\ldots,L.
\]
For $\rho_\ell\geq\bar\rho_\ell$, all the relaxation constraints are
active at an optimal solution, and the lifted variables coincide with
the forward pass of the original network. Hence
\[
G_{\boldsymbol\rho}(x)=R(x).
\]

In particular, in the scalar-parameter case
$\rho_1=\cdots=\rho_L=\rho$, it is sufficient to take
\[
\bar\rho
=
\max_{\ell=1,\ldots,L}
\left\{
\|V\|_\infty
L_\sigma^{L-\ell+1}
\prod_{j=\ell+1}^{L+1}
\|W_j\|_{1\to1}
\right\}.
\]
Therefore,
\[
G_\rho(x)=R(x),
\qquad
\forall \rho\geq\bar\rho.
\]

\paragraph{5. Exactness for monotone convex networks.}
Since $\sigma$ is nondecreasing and $W_\ell\geq0$, the mappings
\[
z_{\ell-1}\mapsto
\sigma(A_\ell x+W_\ell z_{\ell-1}+b_\ell)
\]
are componentwise nondecreasing. Hence, among all feasible lifted
variables, the componentwise smallest ones are obtained recursively by
the forward pass,
\[
z_1=\sigma(\phi_1(x)),\qquad
z_\ell=\sigma(\phi_\ell(x,z_{\ell-1})),
\quad \ell=2,\ldots,L.
\]
Since $W_{L+1}\geq0$, $\sigma$ is nondecreasing, and $V\geq0$, this
choice also minimizes the output term $V^\top u$. Moreover, all penalty
terms are nonnegative and vanish at the forward pass. Therefore,
\[
G_{\boldsymbol{\rho}}(x)=R(x)
\]
for every $\boldsymbol{\rho}\geq0$.

\end{document}